\documentclass[10p]{elsarticle}
\usepackage{amssymb, color}
\usepackage{amsmath}
 \usepackage{pifont}
 \usepackage[utf8]{inputenc}
\usepackage{setspace}
\makeatletter
\def\LaTeX{\leavevmode L\raise.42ex
    \hbox{\kern-.3em\size{\sf@size}{0pt}\selectfont A}\kern-.15em\TeX}
\makeatother

\newcommand{\BibTeX}{{\rm B\kern-.05em{\sc
          i\kern-.025emb}\kern-.08em\TeX}}

\makeatletter
\def\@currentlabel{2.1}\label{e:dispaa}
\def\@currentlabel{2.21}\label{e:dispau}
\def\@currentlabel{2.22}\label{e:dispav}
\def\@currentlabel{2.23}\label{e:dispaw}
\def\@currentlabel{2.24}\label{e:dispax}
\def\theequation{\thesection.\@arabic\c@equation}
\makeatother

\renewcommand{\theequation}{\arabic{section}.\arabic{equation}}

\newcommand{\R}{\mathbb R}

\newcommand{\N}{\mathbb N}
\newcommand{\e}{\epsilon}

\def \D{\Delta}
\def \O{\Omega}

\newtheorem{thm}{Theorem} [section]
\newtheorem{lem}{Lemma} [section]

\newtheorem{cor}{Corollary} [section]

\newtheorem{rem}{Remark}[section]
\newtheorem{taggedtheoremx}{Theorem}
\newenvironment{taggedtheorem}[1]
 {\renewcommand\thetaggedtheoremx{#1}\taggedtheoremx}
 {\endtaggedtheoremx}

\renewcommand{\theequation}{\thesection.\arabic{equation}}
\renewcommand{\thesection}{\arabic{section}}
\renewcommand{\theequation}{\thesection.\arabic{equation}}
\let\ssection=\section\renewcommand{\section}{\setcounter{equation}{0}\ssection}
\def \O{\Omega}
\def \p{\partial}

\def \a{\alpha}
\def \b{\beta}
\def \g{\gamma}
\def \l{\lambda}

\begin{document}
 \begin{frontmatter}
\title{  Existence and nonexistence results of polyharmonic boundary value problems with supercritical growth }
  \author[ah1,ah2,ah3]{ Abdellaziz Harrabi}
\ead{abdellaziz.harrabi@yahoo.fr}
\author[mf1,mf2]{Foued Mtiri}
\ead{mtirifoued@yahoo.fr}
\author[f1,f2]{Wafa Mtaouaa}
\ead{mtawaa.wafa@yahoo.fr}
\address[ah1]{Mathematics Department, Northern Border University, Arar, Saudi Arabia.}
\address[ah2]{Institut sup\'{e}rieur de Math\'{e}matiques Appliqu\'{e}es et de l'Informatique, Kairouan, Tunisia.}
\address[ah3]{ Abdus Salam International Centre for Theoretical Physics, Trieste, Italy.}
\address[mf1]{Mathematics Department, Faculty of Sciences and Arts, King Khalid University, Muhayil Asir, Saudi Arabia.}
\address[mf2]{ ANLIG, UR13ES32, Faculty of Sciences of Tunis, Elmanar University, Tunisia.}
\address[f1]{Mathematics Department, Faculty of Sciences and Arts, King Khalid University, Muhayil Asir, Saudi Arabia.}
\address[f2]{AGTS, LR11ES53, Université de Sfax, Tunisie.}

\begin{abstract}
  We establish some existence results of polyharmonic boundary value problems with supercritical growth.  Our approach is based on truncation argument
 as well as $L^{\infty}$-bounds. Also, by virtue of Pucci-serrin's variational identity \cite{PS}, we provide some nonexistence results.
\end{abstract}

\end{frontmatter}

 \section{Introduction}
Consider the following polyharmonic equations $$(P_{m}):\quad(-\Delta)^m u=g(u)\quad \mbox{in }\;\O, $$
under the Dirichlet boundary conditions
\begin{align}
\label{D}
u = \frac{\p u}{\p \nu} = \ldots = \frac{\p^{m-1} u}{\p\nu^{m-1}} = 0 \quad \mbox{on }\; \p\Omega,
\end{align}
or the Navier boundary conditions
\begin{eqnarray}\label{N}
u = \Delta u = \ldots = \D^{m-1} u = 0\quad \mbox{on }\; \p\Omega,
\end{eqnarray}
where $ N\geq 2m+1$, $\Omega\subset \mathbb{R}^{N}$ is a bounded smooth domain and  $g \in C(\R).$ satisfying:
 \begin{align}
\label{cri}
|g(s)|\leq C(1+|s|^{\frac{N+2m}{N-2m}}), \quad  \forall\,\, s  \in \R.
\end{align}
We define the main $m$-order differential operator $D^m$ by
\begin{equation*}
  D^{m}u = \left\{\begin{array}{llll}\nabla\Delta^{j-1}u & \mbox{ if $r=2j-1$},\\
\D^{j}u & \mbox{ if $r=2j.$ }
\end{array}\right.
\end{equation*}
 The appropriate
 functional spaces of the variational setting to the above  problems are respectively
$$ H_0^m(\Omega):=\left\{v\in
H^m(\Omega);\; \nabla^ju=0\;\mbox{on}\;\partial\Omega,\;\mbox{for}\;j=0, 1, . . ,m-1\right\},$$
and
$$ H_\vartheta^m(\Omega):=\left\{v\in H^m(\Omega);\;\Delta^jv=0\;\mbox{on}\;\partial\Omega,\;\mbox{for}\;
j<\frac{m}{2}\right\}.$$
We denote indifferently $H_\vartheta^m(\Omega)$ and $H_ 0^m(\Omega)$  by $H_m$ which is a Hilbert space equipped with the following
scalar product (see \cite{GGS})
\begin{equation*}
\displaystyle\int_{\O} D^{m}u\cdot D^{m}vdx = \displaystyle \int_{\O} \nabla\Delta^{j-1}u\cdot\nabla\Delta^{j-1}v dx\quad \mbox{if }\;m=2j-1,\\
\quad \quad \displaystyle\int_{\O} \D^{j} u\D^{j}v dx\quad \mbox{if }\;m=2j.
\end{equation*}
Set $G(s)= \displaystyle\int^s_0 g(t)dt$. If $g$ satisfies \eqref{cri}, then the associate energy functional of the boundary value problem $(P_{m})$ is defined by
$$ I(u)=\frac{1}{2}\int_\Omega|D^mu|^2dx-\int_\Omega
    G( u)dx,  \;\;\,\forall\, u \in H_m,$$
and belongs to $C^1(H_m)$. So,  $u\in H_m$ is a weak solution of problem  $(P_{m})$ if and only if $u$ is a critical point of $I$, that is
$$\int_\Omega D^mu\cdot D^mvdx= \displaystyle\int_\Omega
    {g}(u)vdx, \;\; \forall\; v\in H_m.$$
If in addition we assume that
$g\in C^{0,a}_{loc}(\mathbb{R})$,  then any weak solution of the boundary value problem $(P_{m})$  belongs in
$C^{2m}(\overline\O)$ (see the appendix of \cite{HA} for more details)\footnote{ This regularity  allows us to recuperate  the last part of the Navier
boundary conditions \eqref{N}}.  The  Palais-Smale compactness condition plays a
central role in the critical points theory and it is satisfied in the most part of literature  under the following standard assumptions \cite{AR}:
\begin{enumerate}
\item[(AR)] Ambrosetti-Rabinowitz condition:  There exist $\b>2$ and $s_0>0$ such that
$$ g(s)s\geq  \b G(s) > 0  \;\;\mbox{ for all } \, |s|>s_0.$$
\item[(SCP)] Subcritical polynomial growth condition: There exist  $1< p<\frac{N+2m}{N-2m}$ and a positive constant  such that $$|g(s)|\leq C(1+|s|^{p})  \;\;\mbox{ for all } \; s\in\mathbb{R}.$$
\end{enumerate}
These assumptions have been relaxed  into the following large subcritical growth conditions \cite{HA}:
\begin{enumerate}
\item [$(G_1)$:] There exist $C>0$ and $s_{0}>0$
such that $$ C|g(s)|^{\frac{2N}{N+2m}}\leq sg(s)-
 2G(s),\, \;\; \forall  \;\; |s|> s_{0},$$
 \end{enumerate}
and
\begin{enumerate}
\item [$(G_2)$:] $\lim\limits_{s\to \infty}
\dfrac{g(s)}{|s|^{\frac{N+2m}{N-2m}}}=0.$
\end{enumerate}
\smallskip

In this paper, we consider
$$
\left\{
\begin{array}{llll}
(-\Delta)^m u=f(u)+\l |u|^{p-1}u, \; \mbox{ in }\; \Omega;\\
  u \mbox{ satisfies } \eqref{D}; \mbox{ or }  \eqref{N},
\end{array}
\right.
\leqno{(E_{m,\lambda})}
$$
where $\lambda$ is a positive real parameter and
\begin{align}\label{psub}
1< p<\frac{N+2m}{N-2m}.
\end{align}
\smallskip

We examine the effect of the parameter $\l$ to study the existence and the nonexistence of regular solutions of $(E_{m,\l})$. In \cite{hbn},  Brezis and Nirenberg considered the problem $E_{1,\lambda}$ with $f(s)= |s|^{\frac{4}{N-2}}s$. They used a variant of the mountain pass theorem without the Palais-Smale condition to  prove the existence of  positive solution for all $\l>0$ if $N\geq 4,$ and $N=3$ with $3<p<5,$ or $\l$ large enough if $N=3$ and $1<p\leq 3$.

When the variational approach cannot be employed, the question of existence of solutions may be
dealt via topological methods, or bifurcation theory, or truncation technique \cite{AS, FLN, GS,  Haj, HA, HRS, RR,  SO}.  The proof of existence is essentially reduced to deriving $L^\infty$-bounds for a priori solutions or
 solutions of the corresponding perturbed equations.  We shall employ truncation procedure together with minimax argument to obtain existence results of problems $(E_{m,\l}),$  for $\l$ large enough where we  {\bf only } assume that $f \in C(\R)$  satisfies  the following assumption near $0$:
\begin{enumerate}
\item [$(H_0)$]: $ \lim\limits_{s\to 0} \frac{f( s)}{s}=L \in [-\infty, \; \lambda_1)$ and there exist $\nu \in [0,1)$ and $C_1>0$ such that
\begin{align}\label{infty}
| f(s)|\leq C_1 |s|^{1-\nu},\quad \; \forall \;\; |s| \leq 1.
\end{align}
\end{enumerate}
\smallskip

Note that if $L \in (-\infty, \;\lambda_1)$, then \eqref{infty} holds with $\nu=0.$
\medskip

Our main existence result reads as follows:
\begin{thm}\label{th:3} Assume that $p$ verifies \eqref{psub} and $f \in C(\R)$ satisfies $(H_0)$. Then, there exist $\underline{\nu}=\underline{\nu}(\N,p,m)>0$ and $\overline{\l}=\overline{\l}(\nu_0,\O,N,p,f,m)>0$ such that for all $\nu <\underline{\nu}$ and $\l\geq \overline{\l}$:
\begin{enumerate}
 \item  The problem $(E_{m,\l})$
 admits  a nontrivial solution belonging in $C^{2m-1}(\overline{\O})$ and $\|u\|_{L^\infty(\O)}< 1.$  Moreover, if
  $f \in C^{0,a}\left((-1,1)\right)$ for some $a\in (0,1)$, then $u \in$ $C^{2m}(\overline{\O}).$
 \item  If we assume in addition that $f(s)s>0,\,\ \,\forall \,\  |s|\leq 1 \,\ \, \mbox{and}\,\ \,  s\neq 0$. Then the  problem  $(E_{m,\l})$ with  the  Navier  boundary conditions \eqref{N}, admits two nontrivial regular solutions $u_-<0<u_+$.
\end{enumerate}
\end{thm}
 Consider the following Dirichlet boundary value problem
 \begin{equation} \label{n}
(-\D )^mu=|u|^{q-1}u+\lambda |u|^{p-1}u \; \,\,\mbox{in}\,\,\,\O, \,\,\,\, u  \mbox{ satisfies } \eqref{D},
\end{equation}
 When $\O$ is a star shaped domain and $m=1$, the first relevant nonexistence result follows from Pohozaev's identity for all $q\geq \frac{N+2}{N-2}$ and $\lambda\leq0$ (see \cite{cjl, P}).  For more general higher order Dirichlet boundary value problem, Pucci and Serrin \cite{PS} established a  variational identity  allowing further nonexistence results  . In particular they proved
 \begin{taggedtheorem}{A}
\begin{itemize}\item Let $u\in C^{2m}(\Omega)\cap C^{2m-1}(\overline{\Omega})$ be a solution of the Dirichlet boundary value problem $(P_{m})$. Assume that there exists $q >\frac{N+2m}{N-2m}$ such that $ g(s)s-  (q+1)G(s)\geq 0$ for all $s\in \R$.  Then $u\equiv0.$

\item Moreover, if $u$ is a solution of problem \eqref{n} with $\lambda<0,$ $q\geq\frac{N+2m}{N-2m};$ or $\lambda=0$ and $q>\frac{N+2m}{N-2m}.$  Then $u\equiv0.$
 \end{itemize}
\end{taggedtheorem}
The critical growth $q= \frac{N+2m}{N-2m}$ is more difficult even for the case $\lambda=0,$  which has been achieved if $m=2$ for only positive solutions of \eqref{n} \cite{OS} or nodal solutions in  a ball \cite{GR}.  When $\lambda>0$, $m=1$ and $\O$ is a ball, Brezis and Nirenberg \cite{hbn} proved that nontrivial positive bounded solutions of \eqref{n} only exist in the interval $(\frac{\lambda_{1}}{4}, \lambda_{1})$ for $N=3$ and $(0, \lambda_{1})$ for all  $N\geq4$ (see also \cite{FVH1} concerning radial solutions  with prescribed number of nodes). For the higher order case $m\geq 2$, similar result was stated  in only a so called critical  dimension $2m+1\leq N\leq 4m-1$ \cite{GA, GR1, PS1}. The reader may consult \cite{GGS0, GGS} for further comments including also some nonexistence results involving the Navier boundary value biharmonic problem.
\medskip

In the following,  we assume that
\begin{align}\label{psub2}
1\leq p<\frac{N+2m}{N-2m},
\end{align}
and
\begin{itemize}
\item [$(H_1)$]: There exist $q >\frac{N+2m}{N-2m}$, $ s_0>0$  and $C_0>0$ such that
$$f(s)s-  (q+1)F(s)\geq 0, \quad \mbox{for all} \quad s \in \R \quad\mbox{and}\quad f(s)s\geq C_0|s|^{q+1},\quad\mbox{for all} \quad   |s|\leq s_0.$$
\end{itemize}
We exhibit a $\underline{\l}>0$ such that a necessary condition for a nontrivial solution of the Dirichlet boundary value problem  $(E_{m,\l})$ to exist is $\l> \underline{\l}$. This extends the nonexistence results stated in Theorem {\bf A}.
\begin{thm}\label{th:31}
Let $\O$ be a smooth star shaped domain.  Assume that $p$ verifies \eqref{psub2} and $f \in C(\R)$ satisfies $(H_1)$.  Then, there exists $\underline{\l}=\underline{\l}(p,q,m,s_f,N, \O)>0$ such
that for any $\l < \underline{\l},$  the  problem  $(E_{m,\l})$ with the Dirichlet boundary conditions \eqref{D},  has no nontrivial regular solution.
\end{thm}

\begin{rem}
\begin{itemize}
\item  Assume that $p$ verifies \eqref{psub} and $f \in C(\R)$ satisfies $(H_0)$ with $L\in [0,\l_1),$ and $(H_1)$. Thanks to Theorems \ref{th:3} and \ref{th:31} we can find $\overline{\l}>0$ and $\underline{\l}>0$  such that the Dirichlet boundary value problem $(E_{m,\l})$ admits a nontrivial  solution if $\l > \overline{\l}$ and has no a nontrivial solution if $\l < \underline{\l}.$
\item  The following nonlinearities satisfy $(H_0)$-$(H_1)$ $$f(s):= |s|^{q-1}s\exp(as),\,\  a\geq 0 \quad \mbox{ and } \quad f(s):= Ls\exp( (q+1)s)\quad \mbox{ with } \quad q>\frac{N+2m}{N-2m} \mbox{ and }L\in [0,\l_1).$$
\end{itemize}
\end{rem}

Next, consider the following problem
  $$
\left\{
\begin{array}{llll}
(-\Delta)^m u=\a^{-1} f(\a u)+  |u|^{p-1}u\; \mbox{in }\; \Omega; \\
  u \mbox{ satisfies } \eqref{D}; \mbox{ or }  \eqref{N},
\end{array}
\right.
\leqno{(Z_{m,\a})}
$$
 where $\a \in (0,1)$.  From Theorems \ref{th:3} and \ref{th:31}, we derive the following Corollary.

\begin{cor}\label{co:4}
   Assume that $f \in C(\R)$ and $p$ satisfy respectively $(H_0)$ and  \eqref{psub}.  Then there exist $\underline{\nu}=\underline{\nu}(,N,m,p)>0$ and $\underline{\a}=\underline{\a}(f,N,m,\O)>0$ such that for any $0 \leq \nu\leq\underline{\nu}$ and $0 \leq \a\leq\underline{\a},$ we have
 \begin{enumerate}
 \item The  problem $(Z_{m,\a})$
 admits a nontrivial solution belonging in $C^{2m-1}(\overline{\O})$ and $\|u\|_{L^\infty(\O)}< \a^{-1}.$  Moreover, if
  $f \in C^{0,a}\left((-1,1)\right)$ for some $a\in (0,1)$, then $u \in$ $C^{2m}(\overline{\O}).$
 \item  If we assume in addition that $f(s)s>0,\, \forall |s|\leq 1\,\ \, \mbox{and}\,\ \,  s\neq 0$. Then the Navier boundary value problem  $(Z_{m,\a})$
 admits two nontrivial regular solutions $u_-<0<u_+.$
\item Let $\O$ be a star shaped domain.  Assume that $f \in C(\R)$ and $p$ satisfy respectively $(H_1)$ and  \eqref{psub2}.  Then there exists $\overline{\a}=\overline{\a}(p,q,m,f,N, \O)>0$ such
that for any $\a > \overline{\a},$  the Dirichlet boundary value problem  $(Z_{m,\a})$  has no nontrivial regular solution.
\end{enumerate}
\end{cor}
\medskip

The proof of Corollary \ref{co:4} follows from a simple scaling argument. In fact, set $v =\a^{-1} u$ with $\a=\l^{\frac{- 1}{p-1}}$, it is easy to see that $v$ is a solution of $(Z_{m,\a})$ if and only if $u$ is a solution of $(E_{m,\l})$.
\smallskip

\begin{rem}
\begin{itemize}
\item Point 1 of this Corollary improves the existence result stated in \cite{AS}.
\item Set $f(s)=\frac{ -|s|^{2-\nu}\exp( s)}{s}$. Then $f$ satisfies $(H_0)$ with $L =-\infty$ and Corollary \ref{co:4} implies that the following problem has a nontrivial solution $$(-\Delta)^m u=-\frac{|s|^{2-\nu}\exp(\a s)}{\a^\nu s}+  |u|^{p-1}u\;\;\; \mbox{in }\;\;\; \Omega,$$ for $\nu$ and $\a$ small enough.  Point out that the above equation cannot be seen as a small perturbation of the following subcritical problem $$(-\Delta)^m u=  |u|^{p-1}u\;\;\; \mbox{in }\;\;\; \Omega.$$
\end{itemize}
\end{rem}
\medskip

This paper is organized as follows:  In section 2, we provide the proofs of  Theorems \ref{th:3}. Section 3 is devoted to the proofs of  Theorems \ref{th:31}.
\medskip

\section{ Proof of Theorem \ref{th:3}}
\medskip

In the following $C $ denotes always generic positive constants depending on $(N, \O,p,m,C_1)$ only, which could be changed from one line to another.
\medskip

{\bf Proof of point $1$.} According to assumption $(H_0)$, we can find $0<\e_0 <\inf(1,\l_1)$ and $0<s'_0 <1$ such that
\begin{align}\label{F00}
  F(s)\leq \frac{\lambda_1-\e_0}2 s^2  \quad \; \forall \;\; |s| \leq s'_0,
\end{align}
and
\begin{align}\label{h00}
 | f(s)s|\leq C |s|^{2-\nu}, \quad \mbox{ and } \quad  |F(s)|\leq  C |s|^{2-\nu}  \quad \; \forall \;\; |s| \leq s'_0.
\end{align}
Observe that we may rewrite the Poincar\'{e} inequality as follows
\begin{equation}\label{e:dc}
\frac{1}{2}\int_\O  |D^m u|^2dx-\frac{\lambda_1-\e_0}2 \int_\O u^2dx\geq \frac {\e_0} {2\lambda_1}\int_\O  |D^m u|^2dx ,\quad \forall \;\; u\in H_m.
\end{equation}
Let $\theta\in C^1(\R)$ the cut-off function such that
$\theta (s)= 1\;\;  \mbox{if }\;\;  |s|\leq \frac{s'_0}{2}, \quad
 \theta (s)=0\;\;  \mbox{if }\;\; |s|\geq\ s'_0, \, \mbox{ and  }
 0\leq\theta\leq 1 \;\forall s\in \R.$
For $\a \in (0,1)$, set
$$f_\a(s)=\frac{\theta(\a s)}{\a}f(\a s), \quad F_\a(s)=\int_0^s f_\a(t)ds, \quad g_\a(t)= f_\a(s)+| s|^{p-1}s,$$
and $G_\a(s)=\displaystyle\int_0^s g_\a(t)dt = F_\a(s) +\frac{| s|^{p+1} }{p+1}.$
In view of the definition of $\theta$ and \eqref{F00} (respectively \eqref{h00}), we deduce then
\begin{equation}\label{e:kt}
 0 \leq F_\a(s)=\frac{ F(\a s)}{\a^2}\leq \frac{\lambda_1-\e_0}2 s^2, \;\;\;\; \forall \;\; (\a,s)\in (0 , 1)\times
 \R,
\end{equation}
respectively
\begin{equation}\label{e:kt0}
 |f_\a(s)|\leq C \a^{-\nu} , \;\;|f_\a(s)s|\leq C \a^{-\nu}\quad \mbox{ and }\quad |F_\a(s)|\leq C \a^{-\nu}, \;\;\, \forall \;\; (\a,s)\in (0 , 1)\times
 \R.
 \end{equation}
 Consider the truncated problems
 under the Navier or the Dirichlet boundary conditions
$$
\left\{
\begin{array}{llll}
(-\Delta)^m u=g_\a(u) \;\;\; \mbox{in }\; \Omega;\\
 u \mbox{ satisfies } \eqref{D}; \mbox{ or }  \eqref{N}.
\end{array}
\right.
\leqno{(E_{m,\a})}
$$
The associate  energy functional of problems $(E_{m,\a})$ is
$$ I_\a(u)=\frac{1}{2}\int_\Omega|D^m u|^2dx-\int_\Omega G_\a( u)dx, \quad  u \in H_m.$$
As $g_\a(s)= | s|^{p-1}s $ for all $|s|\geq \frac{s'_0}{\a}$, then $g_\a$ satisfies $(AR).$ Also  \eqref{h00} implies
\begin{align*}
(SCP)_\a:\quad  |g_\a(s)|\leq C\a^{-\nu}(1+|s|^{p}) \quad \mbox{for all} \quad s\in\mathbb{R}.
\end{align*}
 Consequently, $ I_\a\in C^1(H_m)$  and satisfies the Palais-Smale condition.  From \eqref{e:dc} and \eqref{e:kt}, one concludes
$$ I_\a(u)\geq \frac {\e_0} {2\lambda_1}\int_\Omega|D^mu|^2dx -\frac{1 }{p+1}\int_\Omega| u|^{p+1}dx.$$
We apply Sobolev's inequality, there holds
$$ I_\a(u)\geq \frac {\e_0} {2\lambda_1}\int_\Omega|D^mu|^2dx -C\left(\int_\Omega|D^mu|^2dx\right)^\frac{p+1}2.$$
Hence, there exist two constants  $\rho,\,\b>0$ such that $(I_\a)_{\mid \partial B_\rho}\geq \b.$ Fix $w_0 \in H_m$, $w_0\neq 0$ and set $v_\nu=b_0\a^{-\frac{\nu}{p+1}} w_0$, $b_0>0$. As $p+1>2$ and $\a \in (0,1)$, from \eqref{h00}, we can find $b_0>0$ large enough such that
$$ I_\a(v_\nu)\leq \a^{-\nu}\left(\frac{b_0^2}{2} \int_\O |D^m w_0|^2dx-\frac{b_0^{p+1}}{p+1}\int_\O w_0^{p+1}dx+C\right)\leq 0.$$
According to the mountain pass theorem, $I_\a$ admits a nontrivial critical point $u_\a \in H_m$ satisfying $$I_\a(u_\a)= \inf\limits_{\gamma\in\Gamma}\sup\limits_{t\in[0,1]} I_\a(\gamma(t))\quad \mbox{where }\quad\Gamma= \left\{\gamma\in C([0,1], H_m);\;
\gamma(0)=0, \gamma (1)=v_\nu\right\}.$$ Using again \eqref{h00}, we derive
\begin{align}\label{za}
I_\a(u_\a)\leq \sup_{t\in[0,1]} I_\a(tv_\nu)\leq \frac{t^2}{2} \int_\O |D^m v_\nu|^2dx -\int_\O F_\a(tv_\nu)dx \leq C\a^{-\nu}.
\end{align}
As $$I'_\a(u_\a)u_\a=\displaystyle\int_\Omega |D^mu_\a|^2dx- \displaystyle\int_\Omega {g}(u_\a)u_\a dx=0,$$ we obtain
\begin{align*}
\begin{split}
& \;I_\a(u_\a)=I_\a(u_\a)- \frac{1}{p+1}I_\a'(u_\a)u_\a \\
&=\left(\frac{1}{2}-\frac{1}{p+1}\right)\int_\O |D^mu_\a|^2dx-\int_\O F_\a(u_\a)dx
+ \frac1{p+1}\int_\O u_\a f_a(u_\a)dx\leq C\a^{-\nu}.
\end{split}
\end{align*}
According to \eqref{za} and  \eqref{h00}, there holds

\begin{align}\label{eg}
\int_\O |D^mu_\a|^2dx \leq C\a^{-\nu}.
\end{align}
Now, we need  the following Lemma based on a  boot-strap argument:
\begin{lem}\label{L}
  $u_\a \in C^{2m-1}(\overline{\O})$ and there exist $C>0$ and  $\g>0$ depending respectively only on $(\O,N,m,p,C)$  and  $(N,p,m)$ such that
 \begin{align}\label{main}
\|u_\a\|_{C^{2m-1,\a'}(\O)} \leq  C\a^{-\nu\g}.
\end{align}
\end{lem}
Precisely if $\frac{2N}{p(N-2m)}<\frac N{2m},$  then $\g=\frac {2mp^3}{N(p-1)}\left( \frac {2m}{N(p-1)}-\frac {(N-2m}{2N}
    \right)^{-1}.$
\medskip

\noindent{\bf Proof.}  Sobolev's inequality gives $$ \|u_\a\|_{L^\frac{2N}{N-2m}(\O)}\leq A \|D^m u_\a\|^2_{L^2(\O)},$$ which combined with  $(SCP)_\a$ imply
 \begin{align}\label{Les1}
  \|g_\a(u_\a)\|_{L^{q_1}(\O)}\leq C \a^{-\nu}\left(1+\|D^m u_\a\|^{p}_{L^2(\O)}\right), \quad \mbox{where} \quad q_1=\frac{2N}{p(N-2m)}.
\end{align}
We invoke now Rellich-Kondrachov's theorem and $L^q$-$W^{2m,q}$ regularity \cite{GGS} \footnote{ Note that the $L^q$-$W^{2m,q}$ regularity is also valid under the Navier boundary conditions (see \cite{Haj}).}.  Precisely, {\bf if }$g_\a(u_\a) \in L^q(\O)$ for some $q>1$, we have
\begin{align}\label{3.5}
\|u_\a \|_{L^{q^*}(\O)}\leq C\|u_\a\|_{W^{2m,q}(\O)} \leq C  \|g_\a(u_\a) \|_{L^q(B_1)},
\end{align}
where
$$ q^*=\frac{ qN}{N-2mq} \quad\mbox{if} \;\; 2mq <N, \quad \mbox{ and } \quad q^*= p\frac{N+1}{2m} \quad\mbox{ if }\;\; q=\frac{N}{2m},$$
and
\begin{align}\label{3.6'}
\|u_\a\|_{C^{2m-1,a'}(\O)} \leq C\|u_\a\|_{W^{2m,q}(\O)}\leq C  \|g_\a(u_\a) \|_{L^q(B_1)} ,\quad \mbox{ if } \;\; 2mq > N.
\end{align}
So, if $q>p$,  $(SCP)_\a$ and \eqref{3.5} imply
 \begin{eqnarray}\label{3.6}
\|g_\a \|_{L^{\frac{q^*}{p}}(B_\frac 12)} \leq C \a^{-\nu}\left(1+ \|g \|^{p}_{L^{q^*}(B_1)}\right), \quad \quad \mbox{ if }\;\; 2mq \leq N.
\end{eqnarray}
Consequently, if $2mq_1> N $ (respectively $2mq_1= N $, the desired estimate \eqref{main} follows from \eqref{Les1} and \eqref{3.6'}  (respectively \eqref{3.6} with $q=q_1$, \eqref{3.6'}-\eqref{3.6} with $q= p\frac{N+1}{2m}$  and \eqref{Les1}).
\smallskip

Let us now give more attention to the difficult case $2mq_1<  N $. Firstly, we have  $$q_1^*=\frac{ q_1N}{N-2mq_1}=\frac{ 2N}{p(N-2m)-4m}>\frac{ 2N}{N-2m}>p\quad\mbox{as}\quad p(N-2m)<N+2m.$$
Set  $q_2=\frac{q_1^*}{p}$, $q_{k+1}=\frac{q_{k}^*}{p}$. We claim that there exists $k_0\in \N^*$ such that
\begin{align}\label{3.8}
 2mq_{k_0+1} >N \quad\mbox{and}\quad  2mq_{k_0} <N.
\end{align}
Indeed, working by contradiction and suppose that $2mq_{k} <N$ for all $k \in \N^*$. Recall that $q^*=\frac{ qN}{N-2mq}\;\; \mbox{if}\;\; 2mq <N $, then $\frac {1}{q_{k+1}}=\frac {p}{q_{k}}-\frac {2mp}{N}$.
From a direct calculation, we obtain
\begin{align}\label{3.11}
  \frac {1}{q_{k+1}}=\frac {p^{k}}{q_{1}}-\frac {2mp}{N}\sum_{ j=0}^{k-1}p^{j}=p^{k}\left(\frac {1}{q_{1}}-\frac {2mp}{N(p-1)}\right)
  +\frac {2mp}{N(p-1)}.
\end{align}
As $p$ verifies \eqref{psub} and  $q_1=\frac{2N}{p(N-2m)}$, then $\frac {1}{q_{1}}-\frac {2mp}{N(p-1)}<0$.  Consequently $\frac {1}{q_{k}}$ converges to $-\infty$, so we reach a contradiction since we assume that $\frac {1}{q_{k}}>\frac {1}{2mN}$. Set   $$\b=\frac {2mp}{N(p-1)}\left( \frac {2mp}{(p-1)}-\frac {1}{q_{1}}
    \right)^{-1} .$$
From \eqref{3.11}, we may see that
 \begin{align}\label{3.12}
 p^{k_0}< \b \quad  \mbox{and}\quad p^{k_0+1}>\b.
 \end{align}
 Therefore, iterating \eqref{3.6} and using \eqref{3.12}, there holds
 \begin{align*}
\|g_\a(u_\a)\|_{L^{q_{k_0+1}}( \O)} \leq  C\a^{-\nu^{p^{k_0}}}\left(1+\|g_\a(u_\a) \|^{p^{k_0}}_{L^{q_1}(\O)} \right)\leq C\a^{-2\nu^{p^{k_0}} } \left(1+\|D^m u _\a\|^{p^{k_0+1}}_{L^2(\O)}\right).
\end{align*}
 From \eqref{3.8}  we have $mq_{k_0+1}> N$ , then the last inequality combined  with \eqref{3.6'}, \eqref{Les1} and \eqref{eg} imply
\begin{align*}
\|u_\a\|_{C^{2m-1,a'}(\O)} \leq   \|g_\a(u_\a) \|_{L^{q_{k_0+1}}(B_1)} \leq  C\a^{-2\nu^{p^{k_0}}} \left(1+\|D^m u _\a\|^{p^{k_0+1}}_{L^2(\O)}\right)\leq  A\a^{-2\nu^{2p^{k_0+1}}}.
\end{align*}
Therefore, $u_\a \in C^{2m-1,a'}((\overline{\O})$ and the desired inequality \eqref{main} follows  from \eqref{3.11}.\qed
\medskip

Set $\underline{\nu}=\g^{-1}$ and take if $\nu < \underline{\nu}$, then from \eqref{main} we can choose  $\a_0 \in (0,1)$ small enough such that
  $$C \a \a^{-\g\nu} < \frac{s'_0}2 \quad\mbox{ so that }\quad|\a u_\a |< \frac{s'_0}{2} \quad
\mbox{ and }\quad g_\a(u_{\a})=\frac1{\a}f(\a u_{\a })+|u_{\a}|^{p-1}u_{\a},\quad\forall\;  0<\a<\a_0.$$
Set $v_{\a}=\a u_{\a},\;$ $\overline{\l} = \a_0^{1-p}$ and $\l = \a^{1-p}$, we obtain
\begin{equation*}
 (-\Delta)^{m} v_{\a} = \a g_\a(u_\a) =
f(v_{\a})+ \l|v_{\a}|^{p-1}v_{\a}.
\end{equation*}
So, $v_\a \in C^{2m-1,\a'}((\overline{\O})$ is a solution of problem $(E_{m,\l})$ with $\lambda = \a^{1-p}$ and $\|v_\a\|_{L^\infty(\O)}\leq 1.$ This means that Theorem \ref{th:3} holds for $\l_p = \a_0^{1-p}.$ Moreover, if $f\in C^{0,a}\left((-1,1)\right)$ then $g_\a(u) \in  C^{0,\inf(\a,\a')}(\overline{\O})$ and from Schauder regularity \cite{GGS},  we derive that $v_\a\in C^{2m,\inf(\a,\a')}(\overline{\O}).$ \qed
\medskip

{\bf Proof of point 2.} Recall first the following Lemma \cite{GGS}:
\begin{lem}\label{ L2}
Let $u \in C^{2m}(\overline{\O})$ be a solution of
\begin{align*}
(-\Delta)^m u =h(x),  \mbox{ in } \O, \,\ \ \mbox{$u$ satisfies \eqref{N}},
\end{align*}
where $h \in C(\overline{\O})$ is a nonnegative function such that $h(x_0)>0$ for some $x_0\in \O$. Then, $u(x) >0$ for all $x \in \O.$
\end{lem}
Set $s_{+}= \max(0,s)$, $s_{-}= \max(0,-s)$  and $$g^+_\a(s)= f_\a(s_+) +s_+^p\quad
\mbox{ respectively }\quad g^-_\a(s)= f_\a(-s_-) -s_-^p.$$
  Following the proof of  point $1$, we can show that there exists $\overline{\l}>0$ such that
\begin{align*}
 (-\Delta)^m u=f(u_+) +\l u_+^p, \ \  \mbox{ in } \O, \,\ \ \mbox{$u$ satisfies \eqref{D}},
\end{align*}
respectively
\begin{align*}
 (-\Delta)^m u=f(-u_-) -\l u_-^p, \ \ \mbox{ in } \O, \,\ \ \mbox{$u$ satisfies \eqref{D}},
\end{align*}
admits a nontrivial solution $u_+\in C^{2m-1}(\overline{\O})$ with (respectively $u_-C^{2m-1}(\overline{\O})$) for all $\l \geq \overline{\l}$.  Since  we assume that $f(s)s >0$ for all $s\neq 0$ then  $g^+_\a\geq 0$ (respectively $-g^-_\a\geq 0$). According to Lemma \ref{ L2},  we deduce that $$u_-<0<u_+ \,\ \ \mbox{ on}\,\ \  \O.$$
The proof of Theorem \ref{th:3} is thereby completed.\qed
\section {Proof of Theorem \ref{th:31}}
\setcounter{equation}{0}
As above, $C $  denotes always generic positive constants depending on $(N, \O,p,q,m,f)$.   According to assumption  $(H_1)$, we get
 \begin{align}\label{e:asm}
\left(1-\frac{2N}{(q+1)(N-2m)}\right)f(s)s \leq f(s)s-  \frac{2N}{N-2m}F(s),\;\;\; \forall \; s \in \R,
\end{align}
Also, we can find a positive constant $C=C(s_0,q)$ such that
\begin{align}\label{e1:asm}
C|s|^{q+1} \leq f(s)s,\;\;\; \forall \; s \in \R.
\end{align}
In fact, as $\leq f(s)s,\; \forall \; |s|\leq s_0,$ we have $F(\pm s_0)>0$. From $(H_{1})$  we have $\left(\frac{F(s) }{s^{q+1}}\right)'\geq 0$ for all $s>0$. Therefore,
$$f(s)s \geq F(s)\geq \frac{ F(s_0)}{s_0^{q+1}}s^{q+1}, \;\,\forall\; s\geq s_0.$$
 Observe that the nonlinearity $-f(-s)$ satisfies also $(H_{1})$. Then inequality \eqref{e1:asm} follows.
\medskip
Let $\O$ be a star shaped smooth bounded domain. Then, there exists $y \in \R^n$ such that $(x-y)\cdot \nu(x)\geq0$, for all $x \in \partial \O$. Consider $u \in C^{2m}(\O) \cap C^{2m-1}(\O)$  a solution of  $(P_{m})$ and $a \in \R$. Applying  Pucci-Serrin's variational identity (see \cite{PS} page 702), we obtain

\begin{align*}
\int_{\Omega}\Bigg[ \left(\frac{N}{2}-a-m\right)|D^{m}u|^{2}+aug(u)-NG(u)\Bigg]dx=\frac{-1}{2}\int_{\partial \Omega}|D^m u|^{2}(\nu\cdot (x-y))ds.
\end{align*}
Fix $a=\frac{N-2m }{2}$, we derive
\begin{align}
\label{ Foufou}
\begin{split}
 \int_{\O}\left( g(u)u -\frac{2N}{N-2m} G(u)\right)dx  \leq 0.
\end{split}
 \end{align}
We multiply $(P_m)$ by $u$, then Poincar\'{e}'s and Sobolev's inequalities imply
\begin{align}\label{sopo}
 C \left(\int_{\O}|u|^{p+1}dx\right)^\frac{2}{p+1} \leq \int_{\O}|D^m u|^2 dx=\int_{\O}g(u)udx.
\end{align}
If $ u \in C^{2m}(\O) \cap C^{2m-1}(\O)$ is a nontrivial solution of problem  $(E_{m,\l})$ under the Dirichlet boundary conditions \eqref{D}, then $u$ is also a solution of $(P_m)$ with $g(s)=f(s)+\l|s|^{p-1}s$. Hence,  \eqref{ Foufou}  and  \eqref{e:asm} imply
\begin{align*}
& \; \left(1-\frac{2N}{(q+1)(N-2m)}\right)\int_{\O}f(u)udx\\
  & \;\leq \int_{\O}\left( f(u)u-\frac{2N}{N-2m }F(u)\right)dx
\leq\l\left(\frac{2N}{(p+1)(N-2m)}-1\right)\int_{\O}|u|^{p+1}dx.
\end{align*}
By \eqref{e1:asm}, we get then
\begin{align*}
& \;C\left(1-\frac{2N}{(q+1)(N-2m)}\right)\int_{\O}|u|^{q+1}dx\\
  & \; \leq \left(1-\frac{2N}{(q+1)(N-2m)}\right)\int_{\O}f(u)udx \leq\l\left(\frac{2N}{(p+1)(N-2m)}-1\right)\int_{\O}|u|^{p+1}dx.
\end{align*}
As $\frac{2N}{(p+1)(N-2m)}-1>0,$ we deduce that $\displaystyle\int_{\O}|u|^{q+1}dx=0$ if $\l\leq 0$. Then $u\equiv 0$.
\smallskip

If $\l >0$,  from H\"older's inequality and \eqref{e:asm},  we obtain
\begin{align}\label{e:an}
   \left(\int_{\O}|u|^{p+1}dx\right)^\frac{q+1}{p+1}\leq C \int_{\O}|u|^{q+1}dx\leq C\int_{\O}f(u)u dx \leq  C\l \int_{\O}|u|^{p+1}dx.
 \end{align}
Therefore,
\begin{align}\label{e:am}
\int_{\O}|u|^{p+1}dx \leq C_1\l^\frac{p+1}{q-p},
\end{align}
where $C_1$ is a positive constant depending on $(N, \O,p,q,m,f)$. In view of  \eqref{sopo} and \eqref{e:an}, we get
\begin{align}\label{klM}
  \left(\int_{\O}|u|^{p+1}dx \right)^\frac{2}{p+1}\leq C\int_{\O}|D^m u|^2 dx\leq C \l\int_{\O}|u|^{p+1}dx.
\end{align}
If $p=1$, we get $\l > \frac 1{C}$ . If $p$ satisfies \eqref{psub},  from  \eqref{klM}, one has
\begin{align*}
 C_2 \l^{-\frac{p+1}{p-1}}  \leq \int_{\O}|u|^{p+1} dx,
\end{align*}
where $C_2$ is a positive constant depending on $(N, \O,p,q,m,f)$.  Combining  the last inequality with \eqref{e:am}, we deduce
$$ \l^{\frac{p+1}{q-p}+\frac{p+1}{p-1}}\geq \frac{C_2}{C_1}.$$
In conclusion, the  problem  $(E_{m,\l})$ under the Dirichlet boundary conditions \eqref{D}, has no nontrivial regular solution if $p$ satisfies \eqref{psub} and $\l < \underline{\l} \mbox{ where } \underline{\l}=\left(\frac{C_2}{C_1}\right)^{\delta}\quad \mbox{with} \quad \delta= \left(\frac{p+1}{q-p}+\frac{p+1}{p-1}\right)^{-1}$ (respectively if $p=1$ and $\l<\underline{\l}$ with $\underline{\l}=  \frac 1{C}$). \qed
\medskip

\section*{Acknowledgment}
The authors extend their appreciation to the Deanship of Scientific Research at King Khalid University, Abha, KSA for funding this work through Research Group under grant number (R.G.P-2 / 121/ 42).

\end{document}